\magnification=\magstep1
\font\T=cmr10 at 19pt
\font\Tt=cmr10 at 14pt
\font\t=cmr10 at 12pt
\def\bigtitle{\T}
\def\mtitle{\Tt}
\def\author{\t}

\footnote{} {This paper was written with the help of the Fonds National
Suisse de la
Recherche Scientifique.} \par

\footnote{} {Mathematics subject classification: 14J17 32S25 57M25}

\footnote{} {Key words: non-isolated singularities, Boundary of
the Milnor fiber, Plumbing graph. Waldhausen decompositions.}

\vskip1in

\centerline {\bigtitle The boundary of the Milnor fiber}

\vskip.1in

\centerline {\bigtitle for some
non-isolated germs of complex surfaces. }

\vskip.5in

\centerline{\author  Fran\c coise Michel, Anne Pichon and
Claude Weber }

\vskip.7in

\noindent {\bf Adresses.}

\vskip.1in

\noindent Fran\c coise Michel / Laboratoire de Math\' ematiques Emile Picard  /
Universit\' e Paul Sabatier / 118 route de Narbonne / F-31062 Toulouse / FRANCE

e-mail: fmichel@picard.ups-tlse.fr

\vskip.1in

\noindent Anne Pichon / Institut de Math\' ematiques de Luminy / UPR 9016
CNRS / Case
907 / 163 avenue de Luminy / F-13288 Marseille Cedex 9 / FRANCE

e-mail: pichon@iml.univ-mrs.fr

\vskip.1in

\noindent Claude Weber / Section de Math\' ematiques / Universit\' e de
Gen\`eve / CP
64 / CH-1211 Gen\`eve 4 / SUISSE

e-mail: Claude.Weber@math.unige.ch

\vskip.7in

\centerline {\mtitle Abstract.}

\vskip.3in

We study the boundary $L_t$ of the Milnor fiber for the
non-isolated singularities in ${\bf C}^3$ with equation $z^m -
g(x,y) = 0$ where $g(x,y)$ is a non-reduced plane curve germ. We
give a complete proof that $L_t$ is a Waldhausen graph manifold
and we provide the tools to construct its plumbing graph. As an
example, we give the plumbing graph associated to the germs
$z^2 - (x^2 - y^3)y^l = 0$ with $l \geq 2$. We prove that the
boundary of the Milnor fiber is a Waldhausen manifold new in
complex geometry, as it cannot be the boundary of a normal
surface singularity.

\vskip.5in

\noindent {\mtitle 1. Introduction.}

\vskip.3in

In [M-P] the authors state with a sketch of proof that the
boundary $L_t$ of the Milnor fiber of a non-isolated surface
singularity in ${\bf C}^3$ is a Waldhausen graph manifold
(non-necessarily "reduziert"). These manifolds are conveniently
described by a plumbing graph. In [M-P-W] we determine the
plumbing graph for the boundary of the Milnor fiber of
Hirzebruch singularities $z^m - x^ky^l = 0$. The present paper
is devoted to the study of germs with equation $z^m - g(x,y) =
0$ where $g(x,y)$ is a non-reduced plane curve germ. For them:

\vskip.1in

\noindent 1) We prove in details that $L_t$ is indeed a
Waldhausem manifold (Section 4). The Waldhausen decomposition
for $L_t$ is obtained by gluing two specific Waldhausen
sub-manifolds along boundary torii: the trunk and the
(non-necessarily connected) vanishing zone.

\noindent 2) We prove that the vanishing zone is in fact a
Seifert manifold and we elucidate its structure (Section 5).

\noindent 3) We show how to obtain the trunk (Section 2) and
how to determine the gluing between the two sub-manifolds
(Section 4).

\vskip.1in

Necessary results about Seifert and Waldhausen manifolds are
recalled in section 3. The dictionary which translates
Waldhausen decompositions into plumbing graphs provided by
[N] can then be used to obtain the canonical plumbing graph
for $L_t$.

\vskip.1in

In section 7, the plumbing graph is  given for the
singularities $z^2 - (x^2 - y^3) y^l = 0$ with $l \geq 2$. We
prove that the boundary of their Milnor fiber are Waldhausen
manifolds new in complex geometry, as they cannot be the
boundary of a normal surface singularity. This fact does not
depend on the orientation on $L_t$.

\vskip.1in

Information about the homology of $L_t$ is given in section 8.
For Hirzebruch singularities we obtain the following result.

\vskip.1in

\noindent {\bf Theorem 8.1.} Let $f(x,y,z) = z^m - x^ky^l = 0$
be the equation of a Hirzebruch singularity. Assume that
$gcd(m,k,l) = 1$, that $1 \leq k < l$ and that $m \geq 2$.
Let $d = gcd (k,l)$ and write $\bar k = k\slash d$ and $\bar l
= l \slash d$. Then $H_1 (L_t, { \bf Z })$ is equal to the
direct sum of a free abelian group of rank $2(m-1)(d-1)$ and a
torsion group. The torsion subgroup is the direct sum of
$(m-1)$ cyclic factors. One of them is of order $m \bar k \bar
l$ and the other $(m-2)$ factors are  of order $\bar k
\bar l$.

\vskip.1in

In section 6, we expound when $L_t$ is a lens space for the
germs under consideration in this paper. The reason why lens
spaces come up is explained at the end of section 2.

\vskip.5in

\noindent {\mtitle 2. Definitions and main results.}

\vskip.3in

We consider germs $f(x,y,z) \in {\bf C}\{x,y,z\}$ such that
$f(0,0,0) = 0$. We deal with germs $f$ such that the dimension
of the singular locus $\Sigma (f)$ is equal to $1$. Hence $f$
is reduced. 

\vskip.1in

 We denote by $B^{2n}_r$  the 2n-ball
with radius $r>0$ centered at the origin of ${\bf C}^n$ and by
$S^{2n-1}_r$  the boundary of $B^{2n}_r$. We set $F_0 =
B^6_{\epsilon} \cap f^{-1}(0)$ and $L_0 = S^5_{\epsilon} \cap
f^{-1}(0)$. According to the theory of Milnor [Mi], extended by 
Burghelea and Verona [B-V]) in the non-isolated case, the homeomorphism classes of the pairs
$(B^6_{\epsilon} , F_0)$ and $(S^5_{\epsilon} , L_0)$ do not
depend on $\epsilon > 0$ if it is sufficiently small. As a
consequence, we shall usually remove $" \epsilon"$ from our
notations.

\vskip.1in

The restriction~~
$f\vert B^6_{\epsilon}~\cap~f^{-1}(B^2_{\eta} -\{ 0 \})
\longrightarrow  (B^2_{\eta} -\{ 0 \})$
is a locally trivial diffe\-rentiable fibration whose
isomorphism class does not depend on $\eta > 0$ provided that
$\eta $ is sufficiently small $(0 < \eta << \epsilon)$. See
Milnor  [Mi] and also Hamm-L\^e [H-L]. Therefore, the
diffeomorphism classes of the manifolds
$F_t = B^6_{\epsilon } \cap f^{-1}(t)$ and $L_t =
S^5_{\epsilon} \cap f^{-1}(t)$
do not depend on $t$ if $0 < \vert t \vert \leq \eta $.
We say that $F_t$ is the Milnor fiber of $f$ and that $L_t$ is
the boundary of the Milnor fiber. $F_t$ is oriented by its
complex structure and $L_t$ is oriented as the boundary of
$F_t$.

\vskip.1in

We denote by $n : \tilde F_0 \rightarrow F_0$ the
normalisation. It follows from the arguments in Durfee [D]
that the boundary $\tilde L_0$ of an algebraic neighborhood of
 $n^{-1}(0)$ is well defined. We shall call
$\tilde L_0$ the boundary of the normalisation.

\vskip.1in

The strategy used to obtain the boundary
of the Milnor fiber for non-isolated singularities is the
following. Let
$\Sigma (f)$ be the singular locus of $f$. By hypothesis
$\Sigma (f)$ is a curve. Let $K_0 = L_0 \cap \Sigma (f)$ be the
link of the singular locus in $L_0$. Let $\tilde K_0 = n^{-1}
(K_0)$ be the pull-back of $K_0$ in $\tilde L_0$. A good
resolution of $\tilde F_0$ provides a Waldhausen decomposition
for $\tilde L_0$ as a union of Seifert manifolds such that
$\tilde K_0$ is a union of Seifert leaves. Let $\tilde M_0$ be
a tubular neighborhood of
$\tilde K_0$ in $\tilde L_0$. The closure $\tilde N_0$ of
$(\tilde L_0 - \tilde M_0)$ is called the trunk of $L_t$.
In 4.6 we define a submanifold $M_t$ of $L_t$ called the
vanishing zone around $K_0$. A slighty less general version of
theorem 4.7 can be easily stated as follows.

\vskip.1in

\noindent {\bf Theorem.} \par
\noindent
1) The closure $N_t$ of $L_t \backslash
M_t$ is diffeomorphic to the trunk $\tilde N_0$. 
\par\noindent
2) The
manifold $M_t$ is a Seifert manifold.

\vskip.1in

The construction (see 4.6) of the vanishing zone is so precise
that it gives rise to a very explicit description of $M_t$. To
each irreducible component $\sigma _i$ of the singular locus of
$f$ corresponds a connected component $M_t (i)$ of $M_t$. An
hyperplane section argument provides a plane curve germ ($z^m -
y^{n_i}$) and an integer $k$. Let $d = gcd (n_i , k)$. In section 5 we
prove the following result.

\vskip.1in

\noindent {\bf Theorem 5.4.} The vanishing zone $M_t (i)$
 is the mapping torus of a diffeomorphism $h : G_t \rightarrow G_t$
 such that :

\noindent 1) $G_t$ is diffeomorphic to the Milnor fiber of the
plane curve germ $z^m - y^{n_i}$.

\noindent 2) The diffeomorphism $h$ is finite of order ${n_i
\slash d}$

\noindent 3) If $d < n_i$, $h$ has exactly
$m$ fixed points and the action of $h$ has order ${n_i \slash
d}$ on all other points.

\noindent 4) Around a fixed point $h$ is a rotation of angle
$-2\pi k \slash n_i$.

\vskip.1in

\noindent {\bf Remark.} When $f$ is not analytically equivalent
to $z^m - g(x,y)$ one can have vanishing zones which are
Waldhausen but not Seifert, or Seifert manifolds of a more
complicated nature.

\vskip.1in

It is stated in [M-P] that $L_t$ is never homeomorphic to
$\tilde L_0$. But the particuliar case when $L_t$ is a lens
space is not treated in [M-P] and is rather delicate. To
produce a complete proof of this statement in a forthcoming
paper, the first two authors need a characterization of the
germs $z^m - g(x,y)$ for which $L_t$ is a lens space. Theorem
6.5 solves the problem.

\vskip.1in

\noindent {\bf Theorem 6.5.} The boundary of the Milnor fiber
of a irreducible germ $f(x,y,z) = z^m - g(x,y)$ is a lens space iff $f$ is
analytically equivalent  to $z^2 - xy^l$.

\vskip.1in

\noindent {\bf Remark.} For our purpose lens spaces are defined as
graph manifolds which are obtained from a plumbing graph which is a
"bamboo" with genus zero vertices.

\vskip.1in

For technical reasons, we use in this paper a polydisc
$B(\alpha) = B^2_{\alpha} \times B^2_{\beta} \times
B^2_{\gamma}$ where
$0 < \alpha \leq \beta \leq \gamma \leq \epsilon$
in place of a standard ball $B^6_{\epsilon}$.

\vskip.1in

\noindent {\bf Definition}. The polydisc $B(\alpha)$ is a
Milnor polydisc for $f$ if:

\vskip.1in

\noindent i) For each $\alpha '$ with $0 < \alpha ' \leq
\alpha$ the pair
$(B(\alpha ') ~,~f^{-1}(0) \cap B(\alpha '))$
is homeomorphic to the pair
$(B^6_{\epsilon} ~,~f^{-1}(0) \cap B^6_{\epsilon})$.

\vskip.1in

\noindent ii) For each $\alpha '$ with $0 < \alpha ' \leq
\alpha$ there exists $\eta$ with $0 < \eta << \alpha '$ such
that:

\vskip.1in

\noindent 1) the restriction of $f$ to
$W(\alpha '~, ~\eta) = B(\alpha ') \cap f^{-1} (B^2_{\eta} -
\{ 0 \})$ is a locally trivial differentiable fibration on
$(B^2_{\eta} - \{ 0 \})$

\noindent 2) this fibration does not depend on $\alpha '$ (when
$0 < \alpha ' \leq \alpha)$ up to isomorphism.

\vskip.5in

\noindent {\mtitle 3. Three-dimensional manifolds.}

\vskip.3in

In this section, we recall some facts pertaining to
3-dimensional manifolds in a setting appropriate to our neeeds.

\vskip.1in

We consider differentiable, compact (usually connected)
3-manifolds $M$ possibly with boundary. When the boundary
$\partial M$ is non-empty, we assume that it is a disjoint
union of torii. Manifolds are oriented. Classifications are
done up to orientation preserving diffeomorphism. In the
situations we meet, $M$ is quite often the boundary of a
complex surface $V$. The complex structure gives rise to an
orientation of $V$ and $M = \partial V$ receives an orientation
via the boundary homomorphism $\partial : H_4 (V mod \partial
V; {\bf Z}) \rightarrow H_3 (\partial V; {\bf Z})$.

\vskip.3in

\noindent {\bf 3.1  Seifert foliations.}

\vskip.2in

In this paper, we only need to consider orientable Seifert
fibrations (to be called Seifert foliations, since we have too
many fibrations present). As our manifolds are oriented and
compact, we may define a Seifert foliation on $M$ as an
orientable foliation by circles. Thanks to a theorem of Epstein
[E], this is equivalent to require that there exists a fixed
point free $S^1$-action on $M$ such that the leaves coincide
with the orbits.

\vskip.1in

An exceptional orbit (leaf) is one such that the isotropy
subgroup is non-trivial. It is a finite cyclic subgroup of
order $\alpha \geq 2$. The slice theorem (and orientability of
$M$) imply that for each exceptional leaf there exist:

\noindent i) a tubular neighborhood which is a union of leaves

\noindent ii) an orientation preserving diffeomorphism of this
neighborhood with the mapping torus of a rotation of order
$\alpha$ on an oriented  2-disc, sending leaves to leaves.

\vskip.1in

A Seifert invariant for an exceptional leaf is defined as
follows. Suppose that the rotation angle on the 2-disc is equal
to ${2 \pi {\beta}^* \over \alpha}$. We need the orientation of
the 2-disc to get the correct sign for the angle. We have
$gcd(\alpha , \beta ^*) = 1$ and  we choose
$\beta ^*$ such that $0 < \beta ^* < \alpha$. Now let $\beta$
be any integer
such that $\beta \beta ^* \equiv 1 \pmod {\alpha}$. The pair
$(\alpha , \beta)$ is a Seifert invariant of the exceptional
leaf. See [Mo] pages 135 to 140. The choice of a $\beta
$ in its residue class$\pmod \alpha$ is related to the choice
of a section of the foliation near the exceptional leaf.

\vskip.1in

Let $r \geq 0$ be the number of boundary components of $M$. The
space of leaves is a compact connected orientable surface of
genus $g \geq 0$ with $r$ boundary components.

\vskip.1in

Suppose now that sections of the foliation are chosen on each
boundary component of $M$ and that they are kept fixed during
the following discussion. We then choose a $\beta$ for each
exceptional leaf. Once these choices have been made, the Euler
number $e \in {\bf Z}$ is defined. See [Mo] for details.
Essentially it is the obstruction to extend the section already
defined on some part of the orbit space. The integer $e$
depends on the choice of the $\beta$ 's, but the rational number
$e_0 = e - \sum {\beta _i \over \alpha _i}$ does not. Of
course, if $r > 0$ the numbers $e$ and $e_0$ still depend on
the choice of a section on the boundary of $M$.

\vskip.3in

\noindent {\bf 3.2  Waldhausen manifolds and plumbings graphs.}

\vskip.2in

The manifolds $\tilde L_0$ and $L_t$ we study in this paper are
graphed manifolds in Waldhausen's sense [W]. They will appear in the
following dress.

\vskip.1in

A finite decomposition $M = \cup M_i$ of a 3-manifold $M$  is
Waldhausen if:

\noindent 1) Each $M_i$ is a Seifert manifold

\noindent 2) If $i \neq j$ the intersection $M_i \cap M_j$
 is either empty or  equal to a union of common boundary
components.

\vskip.1in

A manifold is Waldhausen if it admits a Waldhausen
decomposition. It is best described by a plumbing graph. To
begin with, we consider oriented 3-manifolds which are circle
bundles over a closed surface (we only need here to consider
these). Such a bundle is characterised by its Euler number and
the genus of the base space. Two bundles may be glued together
by an operation called plumbing. See [N] for details.

\vskip.1in

A 3-manifold constructed by plumbing is represented by a graph.
The vertices represent the bundles. They carry two integral
weights: the genus $g$ of the base space and the Euler number
$e$. An edge represents a plumbing operation. The dual graph of
a good resolution  for a normal surface singularity is also
weighted like this. If understood as a plumbing graph, it
describes the boundary of a semi-algebraic neighbourhood of the exceptional locus. See [N] for details.
In [N] Neumann assigns a canonical plumbing graph to each
Waldhausen manifold. Particularly useful are the bamboo o--o--
... --o with genus zero vertices  and Euler numbers $e \leq -2$ for
a lens space (see p.327 thm. 6.1) and the star-shaped tree for the
other Seifert spaces (see p.327 cor.5.7 ).

\vskip.3in

\noindent {\bf 3.3  Mapping torii.}

\vskip.2in

\noindent {\bf Definition.} Let $G$ be an oriented
differentiable surface and let $h: G \rightarrow G$
be an orientation preserving diffeomorphism. The mapping torus
$T(h)$ of $h$ is the quotient of the product $G \times
{\bf R}$ by the equivalence relation $(x,t+1) \sim (h(x), t)$.
The manifold $G \times {\bf R}$ is oriented by the product
orientation, $\bf R$ being equipped with the usual one. This
orientation projects down to $T(h)$.

\vskip.1in

The mapping torus $T(h)$ fibers over the circle $S^1$ with fiber $G$
and
$h$ is "the" monodromy (well defined up to isotopy) of this
fibration. Suppose now that $h$ is of finite order. Then the
lines $\{ x \} \times {\bf R}$ in $G \times {\bf R}$
project onto circles in $T(h)$. The mapping torus thus receives
a foliation in circles.  The following
important property of $e_0$ is useful for computations.
For a proof see [P].

\vskip.1in

\noindent {\bf Theorem.} Let $h$ be an orientation preserving
diffeomorphism acting on a closed surface. Then the rational
Euler number $e_0$ of the Seifert foliation on the mapping
torus $T(h)$ vanishes.

\vskip.3in

\noindent {\bf 3.4 Comments.}

\vskip.2in

i) The plumbing graph for $L_t$ can be obtained as follows. The
plumbing graph for the trunk is part of the plumbing graph for the
normalised surface. From the mapping torus of the vertical monodromy, we
obtain the Seifert-Waldhausen invariants of the vanishing zone by
the dictionary given in [P]. Then [N] gives the plumbing graph for
the vanishing zone. The pasting of two Seifert pieces along a common
boundary component is represented in the plumbing graph by a bamboo
having vertices with $g = 0$.

\vskip.1in

ii) Neumann proves in [N] that the boundary of a normal
surface singularity cannot be a non-trivial connected sum.
However, the boundary $L_t$ of the Milnor fiber can be a
non-trivial connected sum, as proved in [M-P]. In this case,
the canonical plumbing graph is non-connected.

\vskip.1in

iii) Usually when lens spaces $L(n,q)$ are considered it is
implicitely assumed  that $n \geq 2$. In this paper we shall call
generalised lens space an oriented 3-manifold which is
orientation preserving diffeomorphic to $L(n,q)$ or $S^3$ or
$S^1 \times S^2$. They are exactly the 3-manifolds which admit a
genus one Heegaard decomposition. A beautiful result of F.
Bonahon [B] says that the Heegaard decomposition is unique up
to isotopy.

\vskip.1in

iv) A manifold which has two Seifert structures (one of them
being non-orientable) is a frequent pebble in the shoe. Let $h$
be "the" orientation-preserving involution of the annulus $S^1
\times [0,1]$. The mapping torus of $h$ is a Seifert manifold
which has two exceptional leaves with $\alpha = 2$. This is the
Seifert structure that Waldhausen calls $Q$. See [W]. We shall not
meet the other Seifert structure.

\vskip.5in

\noindent {\mtitle 4. From the boundary of the normalisation to
the boundary }

\vskip.1in

 \centerline {\mtitle of the Milnor fiber.}

\vskip.3in

Let $g \in {\bf C}\{ x,y \}$ be non-reduced and such that
$g(0,0) = 0$. Let $\prod_{i=1}^l g_i^{n_i}$ be the factorisation
of $g$ into a product of irrreducible factors with $g_i$ prime
to $g_j$ if $i \neq j$. We choose the indices in such a way
that $n_i > 1$ iff $i \leq i_0$ for some $i_0$ with $1 \leq i_0
\leq l$. We choose the coordinate axis such that $x$ is prime
to $g$ and to ${\partial g \over \partial y}$.

\vskip.1in

Now let $f(x,y,z) = z^m - g(x,y)$ and let $\Gamma = \{
{\partial g \over \partial y} = 0 \} \cap \{ f = 0 \}$.  The
singular locus  $\Sigma (f)$ of $f$
 is the intersection of $\{ z = 0 \}$ with $\{ g'(x,y)
= 0 \}$ where $g'(x,y) = \prod_{i=1}^{i_0} g_i^{n_i} $.

\vskip.1in

(4.1) Now let $B( \alpha )$ be a Milnor polyball as defined at
the end of section 2. Let $S$ be the boundary of $B( \alpha )$
and let $S( \alpha ) = S_{\alpha}^1 \times IntB_{\beta}^2 \times
IntB_{\gamma}^2$. We choose $0 < \alpha \leq \beta \leq \gamma
\leq \epsilon$ such that:

\vskip.1in

\noindent 1. $\Gamma \cap S \subset S( \alpha )$ and
$(\{ g = 0\} \cap \{ z = 0 \} \cap S) \subset S( \alpha )$

\vskip.1in

\noindent 2. $L_0 = (\{ f = 0 \} \cap S) \subset \{ \vert z
\vert < \gamma \}$

\vskip.1in

\noindent 3. In $B( \alpha )$ the curve $\Gamma$ intersects
transversally the hyperplanes $H_a = \{ x = a \}$ for all
$a$ with $0 < \vert a \vert \leq \alpha$.

\vskip.1in

Let $F_0 = f^{-1}(0) \cap B( \alpha )$. Then  $L_0 = S \cap
F_0$ is the boundary of $F_0$. The link $K_0$ of the singular
locus $\Sigma (f)$ of $f$ is by definition $K_0 = \Sigma (f)
\cap L_0$.

\vskip.1in

Now let $n: \tilde F_0 \rightarrow F_0$ be the normalisation of
$F_0$. We have seen in section 2 that $\tilde L_0 = n^{-1}(L_0)$
can be identified with the boundary of the normalisation.
Finally let $\tilde K_0 = n^{-1}(K_0)$ be the pull-back of $K_0$
by the normalisation.

\vskip.1in

\noindent {\bf Remark 4.2} The resolution theory implies that
there exists a decomposition of $\tilde L_0$ as a union of
Seifert manifolds such that $\tilde K_0$ is a union of Seifert
leaves.

\vskip.1in

Let $\varphi : {\bf C}^3 \rightarrow { \bf C}^2 $ be the
projection defined by $\varphi (x,y,z) = (x,z)$. For a small
$\theta$ with $0 < \theta << \alpha$ we denote by $M_0$ the
union of the connected components of $\varphi ^{-1}
(S_{\alpha}^1 \times B_{\theta}^2)$ which meet $K_0$.

\vskip.1in

\noindent {\bf Proposition 4.3} There exists a sufficiently
small
$\theta$ such that:

\noindent 1) $M_0 \subset S( \alpha )$

\noindent 2) $M_0 \cap \{ z = 0 \} = K_0$

\noindent 3) $n^{-1}(M_0) = \tilde M_0$ is a tubular
neighborhood of $\tilde K_0$ in $\tilde L_0 $. Moreover $\tilde
K_0$ is the ramification locus of $\varphi \circ n$ restricted
to $\tilde M_0$.

\vskip.1in

\noindent {\bf Corollary 4.4} The closure $N_0$ of $(L_0 -
M_0)$ in $L_0$ is a Waldhausen manifold.

\vskip.1in

\noindent {\bf Proof of the Corollary.} The restriction of the
normalisation $n$ to the closure $\tilde N_0$ of $(\tilde L_0 -
\tilde M_0)$ in $\tilde L_0$ is a diffeomorphism onto $N_0$. But
$\tilde N_0$ is a Waldhausen manifold by remark 4.2.

\vskip.1in

\noindent {\bf Proof of Proposition 4.3.} From (4.1) fact 1. we
have $K_0 \subset S( \alpha )$. Then there exists $\theta$ such
that $M_0 \subset S( \alpha )$. We can choose $\theta$
 small enough such that  $L_0 \cap \{ \vert z \vert \leq \theta
\}$ is a tubular neighborhood of $\{ g = 0 \} \cap L_0$ in
$L_0$. This proves 2). The singular locus of $\varphi$
restricted to $F_0$ is the curve $\Gamma$. Let $\Delta =
\varphi ( \Gamma )$. We can choose $\theta$ still smaller in
order that $\Delta \cap (S_{\alpha}^1 \times B_{\theta} ^2) =
S_{\alpha}^1 \times \{ 0 \}$. As $\Sigma (f) = \{ z = 0 \} \cap
\Gamma$ this proves 3).

\vskip.1in

(4.5) From the definition of $B( \alpha )$ given at the end of
section 1, there exists a very small $\eta$ with $0 < \eta <<
\theta < \alpha$ such that $f$ restricted to $W( \alpha , \eta
) = B( \alpha ) \cap f^{-1} (B_{\eta}^2 - \{ 0 \})$ is a
locally trivial fibration on $(B_{\eta}^2 - \{ 0 \})$. When $0
< \vert t \vert \leq \eta$ we say that $F_t = W( \alpha , \eta
) \cap f^{-1} (t)$ is "the" Milnor fiber of $f$ and that $L_t =
F_t \cap S$ is the boundary of the Milnor fiber of $f$.

\vskip.1in

In $S$ we consider $\bar S (\alpha ) = S_{\alpha}^1 \times
B_{\beta}^2 \times IntB_{\gamma}^2$ and $\bar S (\beta ) =
B_{\alpha}^2 \times S_{\beta}^1 \times IntB_{\gamma}^2$. As
$\alpha , \beta , \gamma$ have been chosen such that $L_o =
(f^{-1} (0) \cap S) \subset (\bar S (\alpha ) \cup \bar S
(\beta ))$ (see (4.1) fact 2) there exists $\eta$ with $0 <
\eta << \alpha$ such that $L_t \subset (\bar S (\alpha ) \cup \bar S
(\beta ))$ for all $t$ with $0 \leq \vert t \vert \leq \eta$.

\vskip.1in

(4.6) Let $M( \eta )$ be the union of the connected components
of $S \cap \{ \vert f \vert \leq \eta \} \cap \{ \vert z \vert
\leq \theta \}$ which meet $K_0$.  Let
$N( \eta ) $ be the closure of $(W( \alpha , \eta ) \cap S ) -
M( \eta )$ in $S$. For any $t$ with $0 \leq
\vert t \vert \leq \eta$ let $M_t = L_t \cap M( \eta )  $ and
let $N_t = L_t \cap N( \eta )$ be the closure of $(L_t - M_t)$
in
$L_t$.

\vskip.1in

\noindent {\bf Theorem 4.7} There exists a sufficiently small
$\eta$ such that for any $t$ with $0 < \vert t \vert \leq \eta$
we have

\noindent 1) $M_t \subset S( \alpha )$

\noindent 2) $f$ restricted to $N( \eta )$ is a fibration on
$B_{\eta}^2$ with fiber $N_t$ for $0 \leq \vert t \vert \leq
\eta$

\noindent 3) $M_t$ has a Seifert structure such that the
restriction of $z$ on any Seifert leaf is constant.

\vskip.1in

\noindent {\bf Remark 4.8}. Theorem 4.7 enables us to describe
$L_t$ as the union of the Seifert manifold $M_t$ with the
manifold $N_t$ which is diffeomorphic to  the Waldhausen
submanifold $\tilde N_0$ of $\tilde L_0$ defined in proposition
4.3. Moreover, the intersection $M_t \cap N_t$ is equal to
$\partial
 M_t = \partial N_t$ which is a disjoint union of torii. Hence
we have:

\vskip.1in

\noindent {\bf Corollary 4.9.} $L_t $ is a Waldhausen manifold.

\vskip.1in

As $f$ induces a deformation retraction of $M( \eta )$
onto the link $K_0$ we say that $M_t$ is the {\bf vanishing
zone} around $K_0$. Up to a diffeomorphism, $N_t$ is a common
Waldhausen submanifold of $L_t$, $L_0$ and $\tilde L_0$. This
is why we say that $N_t$ (resp $\tilde N_0$) is the {\bf trunk}
of $L_t$ (resp $\tilde L_0$). In concrete terms, $L_t$ can be
constructed as the union of $\tilde N_0$ with $M_t$ and small
collars attached on the boundaries. These small collars are
defined with the help of the restriction of $f$ on the
boundary of $N( \eta )$.

\vskip.1in

\noindent {\bf Proof of 4.7.} Proposition 4.3 implies that $M_0
= M( \eta ) \cap f^{-1}(0)$ is included in $S( \alpha )$. As
$S( \alpha )$ is open, we may choose $\eta$ sufficiently small
in order that $M( \eta ) \subset S( \alpha )$.Thus point 1) is
proved.

As noticed  in 4.5, for a sufficiently small and for $t$ such
that $0 \leq \vert t \vert \leq \eta$ we have $L_t \subset \bar
S (\alpha ) \cup \bar S ( \beta )$. Let $L( \eta ) = N( \eta ) \cup M( \eta )$. We restrict $\eta $ to have

$$\Big( L( \eta ) \cap \{ z = 0 \} \cap \{ {\partial g \over \partial
y} = 0 \} \cap \{ \vert x \vert = \alpha \} \Big ) \subset K_0$$

$$\Big( L( \eta ) \cap \{ z = 0 \} \cap \{ {\partial g \over \partial
x} = 0 \} \cap \{ \vert y \vert = \beta \} \Big ) \subset K_0$$

\vskip.1in

The above inclusions imply that the restriction of
$f$ to $L( \eta ) - K_0$ is a fibration.  The boundary of $N(
\eta )$ (which is equal to the boundary of $M( \eta )$) is
included in $S( \alpha )$ and in $\{ \vert z \vert = \theta \}$.
In proposition 4.3 we have chosen $\theta$ such that $\partial
N_0 = \partial M_0$ does not meet $\{ {\partial g \over
\partial y} = 0 \}$. Hence, for a sufficiently small $\eta$ the
boundary of $N( \eta )$ does not meet $\{ {\partial g \over
\partial y} = 0 \}$ either. This proves 2).

We consider the projection $\varphi$ defined in 4.3. For $0 <
\vert t \vert \leq \eta$ let us denote by $\varphi _t$ the
restriction of $\varphi$ to $M_t$. The singular locus of
$\varphi _t$ is $M_t \cap \{ g' = 0 \} = M_t \cap \{ z^m = t\}$.
For each $c$ with $0 \leq \vert c \vert \leq \theta$ we have

$$\varphi ^{-1} (S_{\alpha}^1 \times \{ c \}) = M_t \cap \{ z =
c \} $$.

This gives a foliation in circles on $M_t$ with leaves defined
by $M_t \cap \{ z= c \} $. {\bf This ends the proof of theorem
4.7}.

\vskip.5in

\noindent {\mtitle 5. The vertical monodromy.}

\vskip.3in

With the notations of 4.1, the link $K_0$ of the singular locus
of $f$ has $i_0$ connected components. We choose $i$ with $1
\leq i \leq i_0$ and we denote by $K_i$ the component of $K_0$
which corresponds to the irreducible factor $g_i$ of $g$. More
precisely:

$$K_i = \big( S \cap \{ z = 0 \} \cap \{ g_i (x,y) = 0 \}
\big)$$

Let $M(i)$ be the connected component of the vanishing zone $M(
\eta )$ (see 4.6) which contains $K_i$. Let $\pi : M( \eta )
\rightarrow S_{\alpha}^1$ be the projection on the $x$-axis.
Let $M_t(i) = M_t \cap M(i)$. Let $\pi _t$ be $\pi$ restricted
to $M_t (i)$ with $0 < \vert t \vert \leq \eta$.

\vskip.1in

\noindent {\bf Lemma 5.1}. The projection $\pi _t$ is a
fibration. Moreover the Seifert leaves constructed in 4.7 are
transverse to the fibers of $\pi _t$.

\vskip.1in

\noindent {\bf Proof of lemma 5.1.} The equation of the
singular locus of $\pi _t$ is $\{ z = 0 \} \cap \{ {\partial g
\over \partial y} = 0 \}$. This curve does not meet $M_t (i)$
when $t \neq 0$.

\vskip.1in

We now choose $a$ with $\vert a \vert = \alpha$ and $P \in K_i
\cap \{ x = a \}$. Let $U(P)$ be the connected component of
$\pi ^{-1} (a) \cap M(i)$ which contains the point $P$. Let
 $f_P$  denote $f$ restricted to $U(P)$. Then $f_P$ is a plane
curve germ with an isolated singular point at $P$ and $G_t =
U(P) \cap M_t (i)$ is its Milnor fiber.

\vskip.1in

\noindent { \bf Definition 5.2.} The vertical monodromy around
$K_i$ is the first return diffeomorphism $h : G_t \rightarrow
G_t$ along the Seifert leaves of $M_t (i)$.

\vskip.1in

\noindent {\bf Remark 5.3.} Let $(s^r , w(s))$ be a Puiseux
expansion of the branch $g_i (x,y) = 0$. Then $G_t ' = M_t (i)
\cap \pi^{-1} (a)$ has $r$ connected components. There exists a
monodromy $h' : G_t ' \rightarrow G_t '$ for the fibration $\pi
_t$ such that $(h')^r $ is the vertical monodromy $h$.

\vskip.1in

Consider the following decomposition $g = g_i^{n_i}
\cdot g''$ in ${\bf C} \{x,y \}$ with $g''$ prime to $g_i$.
Let $k$ be the intersection multiplicity  at the origin between
$g_i$ and $g''$. Let $d = gcd (n_i , k)$.

\vskip.1in

\noindent {\bf Theorem 5.4.} The vanishing zone $M_t (i)$
around $K_i$ is the mapping torus of $h : G_t \rightarrow G_t$
 and we have:

\noindent 1) $G_t$ is diffeomorphic to the Milnor fiber of the
plane curve germ $z^m - y^{n_i}$.

\noindent 2) The vertical monodromy $h$ is finite of order ${n_i
\slash d}$

\noindent 3) If $d < n_i$ the vertical monodromy $h$ has exactly
$m$ fixed points and the action of $h$ has order ${n_i \slash
d}$ on all other points.

\noindent 4) Around a fixed point $h$ is a rotation of angle
$-2\pi k \slash n_i$.

\vskip.1in

\noindent {\bf Proof of theorem 5.4.} The fact that the
vanishing zone is the mapping torus of $h$ is an immediate
consequence of lemma 5.1 and definition 5.2.

\vskip.1in

We first prove statements 1) to 4) when $g_i (x,y) = y$. In
this case, $G_t$ is the Milnor fiber of $f(a,y,z) = z^m -
y^{n_i}g''(a,y)$ with $g''$ prime to $y$. Hence $f(a,y,z)$ has
at $P = (a, 0,0)$ the topological type of $z^m - y^{n_i}$. Thus
point 1) is proved. A Seifert leaf of $M_t (i)$ is in the
hyperplane $\{ z = c \}$ with $0 \leq \vert c \vert \leq
\theta$. It is parametrised by $x = ae^{iv}$
with $v \in [0, 2\pi]$. Moreover, there exists a unity $u(a)$
in ${\bf C}\{ a \}$ such that $g''(a,y) = a^k u(a) + y( ... )$.
Hence, the intersection points $(a,y,c)$ of $G_t$ with
this Seifert leaf satisfy an equation of the following type:

$$(\star ) ~~~~  y^{n_i} = \big( a^k u(a) + y ( ... ) \big)^{-1}
(c^m - t )$$

If $y \neq 0$ the order of $h$ is equal to the order of a
rotation of angle $-2\pi k \slash n_i$ on the parametrised
leaf. This order is equal to $n_i \slash d$. Moreover $y = 0$ if
and only if $ c^m = t $. Hence, we have exactly $m$ fixed
points for $h$ when $z$ is equal to each m-th root of $t$. The
equation $( \star )$ gives directly the angle of rotation
around the $m$ fixed points.

\vskip.1in

In the general case, we consider the Puiseux expansion $(s^r ,
w(s))$ of $g_i (x,y)$. If we make the substitution of
variables  $x = s^r$ ,  $y' = y - w(s)$  and $f'(s,y',z) =
f(s^r , y' + w(s) , z)$  we are back to the preceeding case
with $f$ replaced by $f'$.

\vskip.5in

\noindent {\mtitle 6. When is the boundary of the Milnor fiber
a lens space ?}

\vskip.3in
In this Section, we assume that $f$ is irreducible. In 4.8 we have described the boundary $L_t$ of the Milnor fiber
by gluing the vanishing zone $M_t$  to the trunk $N_0 =
\tilde N_0$.

\vskip.1in

\noindent {\bf Proposition 6.1.} 1) A connected component of
$M_t$ is never a solid torus.

2) When $m > 2$ a connected component of $M_t$ has $m$
exceptional leaves or has a basis with non-zero genus or both.

\vskip.1in

\noindent {\bf Proof of proposition 6.1.} In theorem 5.4 we
have described a connected component $M_t (i)$  of $M_t$ as the
mapping torus of the vertical monodromy $h$ acting on
a differentiable surface $G_t$ which is diffeomorphic to the
Milnor fiber of the plane curve germ $z^m - y^{n_i}$ with $n_i
\geq 2$. Hence $G_t$ is always connected and never
diffeomorphic to a disc. As a consequence $M_t (i)$ is never a
solid torus.

When $m > 2$ the surface $G_t$ has non-zero genus. Then:

i) If $h$ is the identity, the basis of $M_t (i)$ is $
G_t $ itself which has non-zero genus.

ii) If $h$ is not the identity, we have proved in 5.4 that $h$
has exactly $m$ fixed points and hence $M_t (i)$ has $m$
exceptional leaves.

\vskip.1in

\noindent {\bf Proposition 6.2.} If $L_t$ is a lens space, then
the trunk $N_0$ is a solid torus and $M_t$ is connected with a
connected boundary.

\noindent {Proof of proposition 6.2.} The boundary components
of a Seifert manifold which is not a solid torus are
incompressible. If the trunk were not a solid torus, $L_t$
would contain incompressible torii.

\vskip.1in

\noindent {\bf Remark 6.3.} By construction,  the number of connected
components of $M_t$ is equal to the number of irreducible components
of the singular locus $\Sigma (f)$ of $f$.

\vskip.1in

\noindent {\bf Corollary 6.4.} If $L_t$ is a lens space, then
$\Sigma (f)$ is an irreducible germ of curve at the origin of
${\bf C}^3$.

\vskip.1in

\noindent {\bf Theorem 6.5.} The boundary of the Milnor fiber
of a irreducible $f(x,y,z) = z^m - g(x,y)$ is a lens space iff $f$ is
analytically equivalent to $z^2 - xy^l$.

\vskip.1in

\noindent {\bf Proof of theorem 6.5.} In [M-P-W] section 4 it is
proved that the lens space $L(2l, 1)$ is indeed the boundary of
the Milnor fiber of $z^2 - xy^l$.

\vskip.1in

Conversely, when $L_t$ is a lens space, propositions 6.1 and
6.2  and corollary 6.4 imply that $m = 2$, that $N_0$ is a
solid torus and that $\Sigma (f)$ is irreducible. Hence we
can write $g(x,y) = g_1 (x,y)^l \cdot g''(x,y)$ with $g_1$
irreducible , $l = n_1 \geq 2$, $g''$ being aither reduced and prime to
$g_1$ or a unity.

Let $ \psi : ({\bf C}^3,0) \rightarrow ({\bf C}^2,0 )$ be the
projection defined by $\psi (x,y,z) = (x,y)$. Let $S_1$ be the
boundary of the polydisc $B_1 = B_{\alpha}^2 \times B_{\beta}^2$
with $0 < \alpha \leq \beta$ such that $B_1$ is a Milnor
polydisc for $g$. Let $K_1 = S_1 \cap \{ g_1 = 0 \}$. By
construction $\psi (M_0 )$ is a tubular neighborhood of $K_1$
in $S_1$ and the closure $W$ of its complement in $S_1$ is
$\psi (N_0)$.

Let us consider the Milnor fibration $\rho = g_1 \slash {\vert
g_1 \vert} : W \rightarrow S^1$ for the plane curve germ $g_1$.
Let $G_1$ be the Milnor fiber of this fibration. Then $\rho
\circ \psi : N_0 \rightarrow S^1$ is a fibration with fiber
$G_1 '$ which is a ramified covering of $G_1$ induced by $\psi$.
The ramification values of this covering are $G_1
\cap \{g'' = 0 \}$.
 Hence the cardinality of the set of
ramification values is equal to the intersection multiplicity
$m_0 (g_1 , g'')$ of
$g_1$ and $g''$ at the origin of ${\bf C}^2$.

As $m = 2$ this covering has degree 2. Hence

$$\chi (G_1 ) = 1 - \mu (g_1)$$

$$\chi (G_1 ') = 2(1 - \mu (g_1)) - m_0 (g_1 , g'')$$

As $N_0$ is a solid torus,  $G_1 '$ is a disjoint union of discs. The only solution
for the second equation just above is $\mu (g_1) = 0$ and $m_0
(g_1 , g'') $ either equals $1$ or $0$.
\par
When $g''$ is not a unity, i.e. $m_0
(g_1 , g'') =1$, then we can choose the axis in such a way that $g_1 (x,y) = y$
and $g'' (x,y) = x$. As a consequence we obtain that
$f(x,y,z) = z^2 - xy^l$.
\par
Otherwise, we can choose the second axis in such a way that $g_1 (x,y) = y$. Then, $f(x,y,z) = z^2 - y^l$ with $l \geq 3$ as $f$ is irreducible. Then the vertical monodromy on the identity ono a surface which has non-zero genus. Then the Vanishing-zone is a Seifert manifold whose basis has non-zero genus. In this case, we never get a lens space. 

\vskip.1in

{\bf End of proof of theorem 6.5.}
\vskip.2in
\noindent
{\bf Remark.} The reducible case $z^2 - y^2$ is treated in [M-P]. It is prove that $L_t$ is then diffeomorphic to 
${\bf S}^2 \times {\bf S}^1$.
\vskip.5in

\noindent {\mtitle 7. Examples.}

\vskip.3in

In this section we apply the method presented above to the
singularities with equation $z^2 - (x^2 - y^3) y^l = 0$ ~~
($l \geq 2$). The  ingredients necessary to get the Waldhausen
structure are stated in proposition 7.1 for $l$ odd and
proposition 7.2 for $l$ even. The proof of these propositions
is immediate from the theorems proved in section 4 and 5. The
invariants of the Waldhausen structure can then be computed
using the classical results recalled in section 3. From these
and from [N] we can get the canonical plumbing graph.

\vskip.1in

\noindent {\bf Proposition 7.1.} Suppose that $l$ is odd and
write $l = 2 \bar l + 1$ ~~$( \bar l \geq 1 )$. Then:

\noindent 1. The trunk is the Waldhausen manifold $Q$.

\noindent 2. The vanishing zone is connected with one boundary
component. More precisely, it is the mapping torus of an
orientation preserving  diffeomorphism $h$ of order $l$ acting
on the Milnor fiber of the plane curve singularity $z^2 - y^l =
0$. It has two fixed points. The rotation angle at the fixed
points is equal to ${- 2 \over l} 2 \pi$. On the complement of
the fixed points the diffeomorphism $h$ induces a free action
of a cyclic group of order $l$.

\noindent 3. The Waldhausen $(\alpha , \beta)$ for the gluing
between the trunk and the vanishing zone is equal to $(l + 3,
1)$.

\vskip.1in

\noindent {\bf Proposition 7.2.} Suppose that $l$ is even and
write $l = 2 \bar l$ ~~$( \bar l \geq 1 )$. Then:

\noindent 1. The trunk is a thickened torus $S^1 \times S^1
\times [0,1]$.

\noindent 2. The vanishing zone is connected and has two
boundary components. More precisely, it is the mapping torus of
an orientation preserving diffeomorphism $h$ of order
$\bar l$ acting on the Milnor fiber of the plane curve
singularity $z^2 - y^l = 0$. Each boundary component of the
fiber is invariant under $h$. The diffeomorphism $h$ has two
fixed points. The rotation angle at each fixed point is
${- 2 \over l} 2 \pi$. On the complement of
the fixed points the diffeomorphism $h$ induces a free action
of a cyclic group of order $l$.

\noindent 3. The Waldhausen $(\alpha , \beta )$ for the gluing
of the two boundary components of the vanishing zone through
the thickened torus is equal to $(l + 3, 1)$.

\vskip.1in

We now describe the plumbing graphs. We call "bamboo" a graph
with the following shape o--o-- ... --o. The length of a bamboo
is its number of vertices.  All vertices in the
plumbing graphs have genus equal to zero. Most of them have
Euler number equal to $-2$. As a consequence, we only point out
Euler numbers which are different from $-2$.

\vskip.1in

To construct the plumbing graph when $l$ is odd, we start with
a bamboo of length $(l + 4)$. At one extremity, we glue two
bamboos of length one. At the other extremity, we glue two
bamboos of length two. The extremity of these glued bamboos has
Euler number equal to $- \bar l$.

\vskip.1in

To construct the plumbing graph when $l$ is even, we start with
a circuit with ($l + 3$) vertices. At one vertex of the circuit,
we glue two bamboos of length one and Euler number $- \bar l -1$.

\vskip.1in

\noindent {\bf Theorem 7.3.} The boundary $L_t$ of the Milnor
fiber of the non-isolated singularity with equation  $z^2 - (x^2
- y^3)y^l$ ~~$( l \geq 2)$ is not diffeomorphic to the boundary
of a normal surface singularity, whatever the orientation on
$L_t$ may be.

\vskip.1in

{\bf Proof of theorem 7.3.} If it were, the quadratic form
associated to the plumbing graph would be definite (negative
definite if $L_t$ is oriented as the boundary of a resulution. See
[H-N-K]). But the graph contains a full subgraph which has an
indefinite quadratic form: the circuit when $l$ is even and the
maximal full subgraph with $-2$ Euler numbers when $l$ is odd. See
[H].

\vskip.5in

\noindent {\mtitle 8. The homology of the boundary of the
Milnor fiber.}

\vskip.3in

\noindent {\bf Theorem 8.1.} Let $f(x,y,z) = z^m - x^ky^l = 0$
be the equation of a Hirzebruch singularity. Assume that
$gcd(m,k,l) = 1$, that $1 \leq k < l$ and that $m \geq 2$.
Let $d = gcd (k,l)$ and write $\bar k = k\slash d$ and $\bar l
= l \slash d$. Then $H_1 (L_t, { \bf Z })$ is equal to the
direct sum of a free abelian group of rank $2(m-1)(d-1)$ and a
torsion group. The torsion subgroup is the direct sum of
$(m-1)$ cyclic factors. One of them is of order $m \bar k \bar
l$ and the other $(m-2)$ factors are  of order $\bar k
\bar l$.

\vskip.1in

The proof is a consequence of the description we give for
$L_t $ in [M-P-W]. The main ingredient is the determination of
the monodromy ${\bf Z} [t, t^{-1}]$ module associated to the
vanishing zone. As we proved in [M-P-W] that $L_t$ is in fact a
Seifert manifold, one can check that the result fits with
[B-L-P-Z].

\vskip.1in

\noindent {\bf Theorem 8.2.} When $l$ is odd, the group
$H_1 (L_t, {\bf Z })$ for the singularity $z^2 - (x^2 -y^3)
y^l$ is  cyclic of order $4l$. When $l$ is even it is the
direct sum of the integers $\bf Z$ and a torsion group of order
$l(l + 3 )$.

\vskip.5in

\noindent {\bf References.}

\vskip.2in

\noindent [B] F. Bonahon: "Diff\'eotopies des espaces
lenticulaires". Topology 22(1983) p.305-314.

\vskip.1in

\noindent [B-L-P-Z] J. Bryden, T. Lawson, B. Pigott and P.
Zvengrowski: "The integral homology of orientable Seifert
manifolds". Topology and its Appl. 127(2003) p.259-275.

\vskip.1in

\noindent [B-V] D. Burghelea and A Verona: "Local homological
properties of analytic sets".

\noindent Manuscripta Math. 7(1972) p.55-66.

\vskip.1in

\noindent [D] A Durfee: "Neighborhoods of algebraic sets".
Trans. Amer. Math. Soc. 276(1983) p.517-530.

\vskip.1in

\noindent [E] D.B.A. Epstein: "Periodic flows on
three-manifolds". Ann. of Math. 95(1972) p.66-82.

\vskip.1in

\noindent [H] F. Hirzebruch: "\"Uber Singularit\" aten komplexer
Fl\"achen". Rend. Math. Appl. V Ser. 25(1966) p.213-232.

\vskip.1in

\noindent [H-N-K] F. Hirzebruch, W. D.  Neumann and S. S. Koh:
"Differentiable manifolds and quadratic forms". Math. Lecturs Notes
vol4 Dekker (New-York) 1972.

\noindent [H-L] H. Hamm and D. T.  L\^e : "Un th\'eor\`eme de
Zariski du type de Lefschetz". Ann. Sci. \'Ecole Norm. Sup.
6(1973) p.317-355.

\vskip.1in

\noindent [J-N] M. Jankins and W. Neumann: "Lectures on Seifert
Manifolds". Brandeis Lecture Notes 2 (March 1983) 83+27p.

\vskip.1in

\noindent [M-P] F. Michel and A. Pichon: "On the boundary of
the Milnor fibre of nonisolated singularities". Int. Math. Res.
Notes 43(2003) p.2305-2311.

\vskip.1in

\noindent [M-P-W] F. Michel, A. Pichon and C. Weber: "The
boundary of the Milnor fiber of Hirzebruch surface
singularities" To be published in the proceedings of the Luminy
meeting of February 2005.

\vskip.1in

\noindent [Mi] J. Milnor, "Singular Points of Complex Hypersurfaces".
Annals of Mathematical Studies 61  Princeton Univ. Press (1968)
122p.

\vskip.1in

\noindent [Mo] J.-M. Montesinos: "Classical Tesselations and
Three-Manifolds". Universitext Series Springer Berlin (1987).

\vskip.1in

\noindent [N] W. Neumann: "A calculus for plumbing applied to
the topology of complex surface singularities and degenerating
complex curves". Trans. Amer. Math. Soc. 268(1981) p.299-344.

\vskip.1in

\noindent [P] A. Pichon: "Fibrations sur le cercle et surfaces
complexes". Ann. Inst. Fourier 51(2001) p.337-374.

\vskip.1in

\noindent [W] F. Waldhausen: "Uber eine Klasse von 3-dimensionalen
Mannigfaltigkeiten". Invent.Math. 3(1967) p.308-333 and 4(1967)
p.87-117.

\vskip.5in

Gen\`eve, le 28 avril 2006.

\end